\documentclass[a4paper,12pt]{article}

\usepackage{color}
\newcommand{\red}[1]{{\colorbox{red}{{#1}}}}

\usepackage{makeidx}

\usepackage{latexsym}
\usepackage{amscd}
\usepackage{amsmath}
\usepackage{amssymb}

\usepackage{amsthm}
\usepackage{float}
\usepackage{graphicx}
\usepackage{psfrag}

\theoremstyle{plain}
\newtheorem{theorem}{Theorem}

\newtheorem{lemma}[theorem]{Lemma}

\newtheorem{proposition}[theorem]{Proposition}

\theoremstyle{definition}
\newtheorem{definition}[theorem]{Definition}
\newtheorem{remark}[theorem]{Remark}
\newtheorem{question}[theorem]{Question}

\newdimen\argwidth
\def\db[#1\db]{%
 \setbox0=\hbox{$#1$}\argwidth=\wd0
 \setbox0=\hbox{$\left[\box0\right]$}
  \advance\argwidth by -\wd0
 \left[\kern.3\argwidth\box0 \kern.3\argwidth\right]}

\newcommand{\Ker}{\operatorname{Ker}}
\newcommand{\Image}{\operatorname{Im}}
\newcommand{\Coker}{\operatorname{Coker}}
\newcommand{\Coim}{\operatorname{Coim}}

\newcommand{\id}{\ensuremath{\mathop{\mathrm{id}}}}
\newcommand{\Sym}{\operatorname{Sym}}

\newcommand{\Spec}{\operatorname{Spec}}

\newcommand{\xto}{\xrightarrow}
\newcommand{\lto}{\longrightarrow}
\newcommand{\hto}{\hookrightarrow}
\newcommand{\simto}{\xrightarrow{\sim}}

\newcommand{\bC}{\ensuremath{\mathbb{C}}}

\newcommand{\bN}{\ensuremath{\mathbb{N}}}

\newcommand{\bP}{\ensuremath{\mathbb{P}}}

\newcommand{\bZ}{\ensuremath{\mathbb{Z}}}

\newcommand{\scE}{\ensuremath{\mathcal{E}}}
\newcommand{\scF}{\ensuremath{\mathcal{F}}}

\newcommand{\scL}{\ensuremath{\mathcal{L}}}

\newcommand{\scN}{\ensuremath{\mathcal{N}}}
\newcommand{\scO}{\ensuremath{\mathcal{O}}}
\newcommand{\scP}{\ensuremath{\mathcal{P}}}

\newcommand{\scT}{\ensuremath{\mathcal{T}}}

\newcommand{\fraka}{\ensuremath{\mathfrak{a}}}

\newcommand{\logf}{\scT_{\bP^\ell}(-\log f)}
\newcommand{\He}{\operatorname{He}}
\newcommand{\bCx}{\bC^{\times}}

\newcommand{\Der}{\mathrm{Der}}
\newcommand{\depth}{\operatorname{depth}}

\title{Logarithmic vector fields along smooth plane cubic curves}
\author{Kazushi Ueda and Masahiko Yoshinaga}
\date{}
\begin{document}

\maketitle

\begin{abstract}
We study the sheaves of logarithmic vector fields
along smooth cubic curves in the projective plane,
and prove a Torelli-type theorem
in the sense of Dolgachev--Kapranov \cite{Dolgachev-Kapranov}
for those with non-vanishing $j$-invariants.
\end{abstract}

\section{Introduction}

K. Saito \cite{Saito_TLDFLVF} introduced the notion of
the sheaf of logarithmic vector fields along a divisor
and proved that it is always reflexive.
A divisor $D$ in a variety $S$ is said to be {\em free}
if the sheaf of logarithmic vector field along $D$
is a free $\scO_S$-module.
He proved that the discriminant in the parameter space
of the semi-universal deformation
of an isolated hypersurface singularity
is always free.

When the ambient space is the projective space $\bP^\ell$,
an $\scO_{\bP^\ell}$-module is said to be free
if it is the direct sum $\bigoplus_{i} \scO_{\bP^\ell}(a_i)$
of invertible sheaves.
The problem of characterizing free divisors
in projective spaces has attracted much attention,
especially when the divisor is given
as an arrangement of hyperplanes.
See e.g. \cite{Terao_GE}.
If a divisor in $\bP^\ell$ is free,
then the passage from the divisor
to the sheaf of logarithmic vector fields
causes loss of information;
only the sequence $\{ a_i \}_{i=1}^\ell$ of integers
is left,
and it is impossible to reconstruct the divisor
from this finite amount of information.

In the opposite extreme,
Dolgachev and Kapranov \cite{Dolgachev-Kapranov}
asked when the the sheaf $\scT( - \log D)$
contains enough information
to reconstruct $D$.
A divisor $D$ in $\bP^\ell$ is said to be {\em Torelli}
if the isomorphism class of $\scT( - \log D)$
as an $\scO_{\bP^\ell}$-module
determines the divisor $D$.
Their main result is the condition
for an arrangement of sufficiently many hyperplanes
in $\bP^\ell$ to be Torelli.

In this paper,
we discuss the case when $\ell = 2$
and $D$ is a smooth cubic curve.
Our main result asserts that
$D$ is Torelli precisely when the $j$-invariant of $D$ is not zero.
The strategy of our proof is the following:
\begin{enumerate}
\item
 The set of jumping lines of the sheaf of logarithmic vector fields
 along a smooth cubic curve coincides with its Cayleyan curve.
\item
 For a smooth cubic curve with a non-vanishing $j$-invariant,
 the Cayleyan curve determines the original curve
 up to three possibilities.
\item
 The set of ``jumping cubic curves'' fixes this left-over ambiguity
 and the Torelli property holds.
\item
 When the $j$-invariant of $D$ is zero,
 we can construct a family of divisors
 with isomorphic sheaves of logarithmic vector fields along them.
\end{enumerate}

Smooth cubic curves with vanishing $j$-invariants
provide examples of divisors
which are neither free nor Torelli.

{\bf Acknowledgment:}
We thank Igor Dolgachev for a stimulating lecture
in Kyoto in winter 2006
and Akira Ishii for valuable discussions and comments.
K. U. is supported
by Grant-in-Aid for Young Scientists (No.18840029).
M. Y. is supported by JSPS Postdoctoral Fellowship 
for Research Abroad. 

\section{Preliminaries}

\subsection{de Rham--Saito's lemma}

Let $A$ be a Noetherian ring and
$
 M = \bigoplus_{i=1}^n A e_i
$
be a free module over $A$ generated by
$e_1, \ldots, e_n$.
For $\omega_1, \ldots, \omega_r\in M$, put
$$
 \omega_1 \wedge \dots \wedge \omega_r
  = \sum_{1 \le i_1 < \dots <i_r \le n}
      a_{i_1, \dots, i_r} e_{i_1} \wedge \dots \wedge e_{i_r}.
$$
and define $\fraka$ to be the ideal generated by
$a_{i_1, \cdots, i_r}$ for $1 \le r \le n$
and $1 \le i_1 < \dots <i_r \le n$.
We also define as follows:
\begin{eqnarray*}
 Z^p
  &=& \{ \varphi \in \wedge^p M \mid
         \omega_1 \wedge \dots \wedge \omega_r \wedge \varphi = 0 \}, \\
 B^p &=& \sum_{k=1}^r \omega_k \wedge (\wedge^{p-1} M), \\
 H^p &=& Z^p / B^p.
\end{eqnarray*}

\begin{theorem}[de Rham--Saito's lemma
\cite{de_Rham_SLDDF,Saito_GdRL}]
\begin{itemize}
\item[(1)] There exists an integer $\nu\in\bZ_{\geq 0}$
such that $\fraka^\nu H^p=0$ for
$0\leq p\leq n$.
\item[(2)] For $0\leq p< \depth_\fraka A$,
we have $H^p=0$.
\end{itemize}
\end{theorem}

\subsection{Sheaf of logarithmic vector fields}

Let $A=\bC[z_0, \ldots, z_\ell]$ be a polynomial ring
and $\Der_A$ be the module of $\bC$-derivations
of $A$,
which is a free module of rank $\ell+1$;
$$
\Der_A=\sum_{i=0}^\ell A\frac{\partial}{\partial z_i}.
$$

\begin{definition}
\normalfont
For a homogeneous polynomial $f\in A$, we define
\begin{align*}
 D(-\log f)
  &= \{ \delta \in \Der_A \mid \delta f \in (f) \}, \\
 D_0(-\log f)
  &= \{ \delta \in \Der_A \mid \delta f = 0 \}.
\end{align*}
We put
$\deg z_i = 1$ and
$\deg(\partial / \partial z_i) = - 1$ for $i=0, \dots, l$.
The degree $k$ part of $D_0(-\log f)$
will be denoted by $D_0(-\log f)_k$.
\end{definition}

We have the direct sum decomposition
$$
 D(-\log f) = D_0(-\log f) \oplus A \cdot E,
$$
where
$$
 E = \sum_{i=0}^\ell z_i \partial/\partial z_i
$$
is the Euler vector field.
Let $\Omega_A$ be the module of differentials
$$
 \Omega_A^1=\bigoplus_{i=0}^\ell A dz_i,
$$
and $\Omega_A^k$ be its $k$-th exterior power
for $k = 0, \dots, \ell+1$.
We have an isomorphism of $A$-modules
$$
 D_0( - \log f)
  \cong \{ \omega \in \Omega^\ell \mid df \wedge \omega = 0 \}
$$
under the identification
$$
\begin{array}{rcl}
 \Der_A & \simto & \Omega^\ell \\
 \rotatebox{90}{$\in$} & & \rotatebox{90}{$\in$} \\
 \sum_{i=0}^\ell f_i \frac{\partial}{\partial z_i} &
 \longmapsto &
 \sum_{i=0}^\ell (-1)^i f_i dz_0 \wedge \dots
   \wedge \widehat{dz_i} \wedge \dots \wedge dz_\ell.
\end{array}
$$
Let $D \subset \bP^\ell$ be the hypersurface
defined by $f$.
If $D$ is smooth,
then the origin $0 \in \bC^{l+1}$ is
the only zero locus of the Jacobi ideal
$$
 J(f) = \left(
         \frac{\partial f}{\partial z_0}, \dots,
         \frac{\partial f}{\partial z_\ell}
        \right),
$$
and hence we have
$$
 \depth_{J(f)} A = \ell + 1.
$$
Let $H^p$ be the $p$-th cohomology
of the complex
$$
 0
  \lto \Omega_A^0
  \xto{df\wedge} \Omega_A^1
  \xto{df\wedge} \cdots
  \xto{df\wedge} \Omega_A^{\ell}
  \xto{df\wedge} \Omega_A^{\ell+1}
  \lto 0.
$$
If $D$ is smooth,
then we have $H^p=0$ for $p = 0, \dots, \ell$
by de Rham--Saito's lemma.
Since
$$
 D_0(- \log f)
  \cong \Ker \left(
              df \wedge : \Omega^{\ell} \to \Omega^{\ell+1}
             \right),
$$
the sequence
\begin{equation} \label{eq:res_Dlog}
 0
  \lto \Omega_A^0
  \xto{df\wedge} \Omega_A^1
  \xto{df\wedge} \cdots
  \xto{df\wedge} \Omega_A^{\ell-1}
  \xto{df\wedge} D_0(- \log f)
  \lto 0
\end{equation}
gives a free resolution of $D_0(- \log f)$.
\medskip

The Euler sequence
$$
\begin{array}{ccccccccc}
 0 & \lto & \scO & \lto & \scO(1)^{\ell+1}
   & \lto & \scT_{\bP^{\ell}} & \lto & 0, \\
 & & \rotatebox{90}{$\in$} & & \rotatebox{90}{$\in$} & & & & \\
 & & 1 & \mapsto & E & & & &
\end{array}
$$
shows that the sheafification $\logf$ of ${D_0(-\log f)}$
can be considered
as a subsheaf of the tangent sheaf $\scT_{\bP^\ell}$;
$$
 \logf \subset \scT_{\bP^{\ell}}.
$$
It is the sheaf of holomorphic vector fields
tangent to the hypersurface $D$ at smooth points of $D$.
If $D$ is smooth,
we have the short exact sequence
$$
 0 \lto \logf \lto \scT_{\bP^\ell}
   \lto \scN_{D/\bP^\ell} \lto 0,
$$
where $\scN_{D/\bP^\ell}$ is the normal bundle.
We have an isomorphism
$$
 d f|_D : \scN_{D / \bP^{\ell}} \simto \scO_D(d),
$$
where
$$
 d = \deg f.
$$

If $D$ is smooth,
then the sheaf $\scT_{\bP^\ell}(-\log f)$ has the resolution
\begin{equation} \label{eq:res_Tlog}
 0 \to \scO(1-(d-1)\ell) \to \dots
   \to \scO(3-2d)^{\oplus \binom{\ell+1}{\ell-2}}
   \to \scO(2-d)^{\oplus \binom{\ell+1}{\ell-1}}
   \to \logf
   \to 0
\end{equation}
obtained by sheafifying the exact sequence
(\ref{eq:res_Dlog}).
We also have
$$
 \Gamma
  \left(
   \bP^\ell, \logf(k)
  \right)
 = D_0(-\log f)_k
$$
for $k \in \bZ$.

\section{Plane curves}

Now we set $\ell=2$
to focus our attention on plane curves.
Let
$
 f \in \bC[z_0, z_1, z_2]
$
be a homogeneous polynomial of degree $d$ and
$D \subset \bP^2$ be the curve
defined by $f$.
Define $\scF$ as the cokernel of
$
 df \wedge : \scO(3 - 2 d) \to \scO(2 - d)^{\oplus 3}
$
so that we have the exact sequence
\begin{equation} \label{eq:exact}
 0
  \lto \scO(3 - 2 d)
  \xto{df\wedge} \scO(2 - d)^{\oplus 3}
  \lto \scF
  \lto 0.
\end{equation}
The Chern polynomial of $\scF(k)$ is given by
\begin{align*}
 c_t(\scF(k))
  &:= 1 + c_1(\scF(k)) t + c_2(\scF(k)) t^2 \\
  &= c_t(\scO(2 - d + k))^3 c_t(\scO(3 - 2 d + k))^{-1} \\
  &= 1 + (3 - d + 2 k) t + (d^2 - 3 d + 3 + k^2 + (3 - d) k)t^2
\end{align*}
for $k \in \bZ$.
If $D$ is smooth,
then we have
\begin{align*}
 \scF
  &:= \Coker
   (df \wedge : \scO(3 - 2 d) \to \scO(2 - d)^{\oplus 3}) \\
  &\cong \Coim
   (df \wedge : \scO(2 - d)^{\oplus 3} \to \scO(1)^{\oplus 3}) \\
  &\cong \Image
   (df \wedge : \scO(2 - d)^{\oplus 3} \to \scO(1)^{\oplus 3}) \\
  &\cong \Ker
   (df \wedge : \scO(1)^{\oplus 3} \to \scO(d)) \\
  &\cong \logf.
\end{align*}

\begin{lemma} \label{lem:stability}
If $D$ is smooth,
then $\scT_{\bP^2}(-\log f)$ is stable.
\end{lemma}

\begin{proof}
We consider $\scF([(d - 3) / 2])$ instead of $\scT_{\bP^2}(- \log f)$
whose first Chern number is normalized to either
$0$ (when $d$ is odd) or $-1$ (when $d$ is even).
Then $\scF([(d - 3) / 2])$ is stable
if and only if it has no global section.
This follows from the cohomology long exact sequence
associated with the short exact sequence (\ref{eq:exact})
tensored with $\scO_{\bP^2}([(d - 3) / 2])$.
\end{proof}

\section{Smooth cubic curves}

Let $f \in \bC[z_0, z_1, z_2]$ be a homogeneous polynomial
of degree three and $D \subset \bP(V)$ be a cubic curve
defined by $f$,
where $V = \Spec \bC[z_0, z_1, z_2]$.
We assume that $D$ is smooth.

\subsection{Jumping lines}

Let $L$ be a point in the dual projective plane $\bP(V^*)$
defined by a linear form
$
 \alpha = \alpha_0 z_0 + \alpha_1 z_1 + \alpha_2 z_2 \in V^*.
$
We can think of $L$ as a line in $\bP(V)$.
Restricting the short exact sequence (\ref{eq:exact}) to $L$
and taking the cohomology long exact sequence,
we have
$$
0
 \lto H^0(\scF|_L)
 \lto H^1(\scO_L(-3))
 \lto H^1(\scO_L(-1))^3
 \lto H^1(\scF|_L)
 \lto 0.
$$
Since
$$
 H^1(\scO_L(-3)) \cong H^0(\scO_L(1))^* \cong \bC^2
$$
and
$$
 H^1(\scO_L(-1)) \cong H^0(\scO_L(-1))^* = 0,
$$
we have
$$
 \dim H^0(\scF|_L) =2.
$$
Hence $\scF|_L$ is either
$$
 \scF|_L =
  \begin{cases}
   \scO_L \oplus \scO_L & \text{$L$ is generic}, \\
   \scO_L(-1) \oplus \scO_L(1) & \text{$L$ is a jumping line.}   
  \end{cases}
$$
In particular,
$$
 \text{$L$ is a jumping line}
  \Longleftrightarrow H^0(\scF(-1)|_L)\neq 0.
$$
By tensoring $\scO_L(-1)$ with
the short exact sequence (\ref{eq:exact})
and taking the cohomology long exact sequence,
we have
$$
\begin{array}{ccccccccccc}
0  &
\rightarrow  &
H^0(\scF(-1)|_L)  &
\rightarrow  &
H^1(\scO_L(-4))&
\xto{df\wedge} &
H^1(\scO_L(-2)^{\oplus 3})&
\rightarrow  &
H^0(\scF(-1)|_L)&
\rightarrow  &
0\\
&&&&\rotatebox{90}{$\cong$}&&\rotatebox{90}{$\cong$}&&&&\\
&&&& H^0(\scO_L(2))^*
&
&H^0(\scO_L^{\oplus 3})^*.&&&&
\end{array}
$$
Since
$
 H^0(\scO(2)|_L)
  \cong \Sym^2V^*/(z_0 \alpha, z_1 \alpha, z_2 \alpha),
$
the set $S = S(\logf) \subset \bP(V^*)$
of jumping lines
is characterized as follows;
\begin{align}
 L \in S
 & \Longleftrightarrow
  \text{$(df \wedge)^* : H^0(\scO_L^{\oplus 3}
   \to H^0(\scO_L(2)))$ is not an isomorphism}
     \nonumber \\
 & \Longleftrightarrow
  \text{
   $z_0 \alpha, z_1 \alpha, z_2 \alpha,
    \partial_0 f, \partial_1 f, \partial_2 f$
   are linearly dependent in $\Sym^2 V^*$}.
     \label{eq:jumping_line}
\end{align}

\subsection{Cayleyan curves}

Here we prove the following:

\begin{proposition} \label{prop:jumping_line-Cayleyan}
Let $D \subset \bP(V)$ be a smooth cubic curve
defined by a polynomial $f$.
Then the set $S = S(\logf) \subset \bP(V^*)$
of jumping lines of $\logf$
in the dual projective plane $\bP(V^*)$
is the Cayleyan curve of $D$.
\end{proposition}

First we recall the definition of the Cayleyan curve
of a plane cubic curve
following Artebani and Dolgachev \cite{Artebani-Dolgachev}.
The {\em first polar} of a plane curve $D = \{ f = 0 \}$
with respect to a point $q = [a_0 : a_1 : a_2] \in \bP(V)$
is the curve
$
 P_q(D)
  = \{
     a_0 \partial_0 f + a_1 \partial_1 f + a_2 \partial_2 f = 0
    \}$
whose degree is one less than that of $D$.
One can show that when $D$ is a cubic curve,
the Hessian curve
$
 \He(D)
  = \{
     \det \left( (\partial_i \partial_j f )_{i,j=1}^3 \right) = 0
    \} \subset \bP(V)
$
consists of points $q \in \bP(V)$
such that the polar curve $P_q(D)$ decomposes
into the union of two lines.
For $q \in \He(D)$,
let $s_q \in \bP(V)$ be the singular point of $P_q(D)$
and $L_q \in \bP(V^*)$ be the line
connecting $q$ and $s_q$.
It is known that $s_q$ lies on $\He(D)$
and the map
$$
\begin{array}{ccccc}
 s & : & \He(D) & \lto & \He(D) \\
 & & \rotatebox{90}{$\in$} & &  \rotatebox{90}{$\in$} \\
 & & q & \mapsto & s_q
\end{array}
$$
is a fixed-point-free involution on $\He(D)$.
The image of the map
$$
\begin{array}{ccc}
 \He(D) & \lto & \bP(V^*) \\
 \rotatebox{90}{$\in$} & &  \rotatebox{90}{$\in$} \\
 q & \mapsto & L_q
\end{array}
$$
is called the {\em Cayleyan curve} of $D$,
which is known to be the quotient of $\He(D)$
by the involution $s$.
A linear form
$
 \alpha = \alpha_0 z_0 + \alpha_1 z_1 + \alpha_2 z_2 \in V^*
$
represents a point in the Cayleyan curve of $D$ 
if and only if there is a point $[a_0 : a_1 : a_2] \in \bP^2$
such that
$$
 a_0 \partial_0 f + a_1 \partial_1 f + a_1 \partial_1 f
  \in \alpha \cdot V^*.
$$
This is precisely the condition (\ref{eq:jumping_line})
for the line $[\alpha] \in \bP(V^*)$
to be a jumping line of $\logf$.

\subsection{The set of jumping lines and $j$-invariant}

Here we prove the following:

\begin{proposition}
Let $D$ be the smooth cubic curve
defined by a polynomial $f$.
Then the set $S(\logf)$ of jumping lines is singular
if and only if the $j$-invariant of $D$ is zero.
\end{proposition}

\begin{proof}
Choose a coordinate of $V$
so that $f$ is a Hesse cubic
\begin{equation} \label{eq:Hesse}
 f_t(z_0, z_1, z_2)
  = z_0^3 + z_1^3 + z_2^3 - 3 t z_0 z_1 z_2,
\end{equation}
where $t \in \bC \setminus \{ 1, \zeta, \zeta^2 \}$
and $\zeta = \exp[2 \pi \sqrt{-1} / 3]$.
Recall that $D = \{ f_t = 0 \} \subset \bP^2$ is smooth
if and only if $t^3 \ne 1$.
The set $S = S(\logf)$ of jumping lines,
which coincides with the Cayleyan curve of $D$,
is a Hesse cubic
$$
 t (\alpha_0^3 + \alpha_1^3 + \alpha_2^3)
  - (t^3 + 2)  \alpha_0 \alpha_1 \alpha_2 = 0
$$
in the dual projective plane.
It is the union of three lines in general position
if $t=0$ or $(3 t)^3 = (t^3 + 2)^3$.
Since
$$
 (t^3 + 2)^3 - (3 t)^3 = (t^3 - 1)^2 (t^3 + 8)
$$
and the $j$-invariant $j(D)$ of $D$ is given by
$$
 j(D) = \frac{1}{64} t^3 \frac{(t^3 + 8)^3}{(t^3 - 1)^3},
$$
the Cayleyan curve of $D$ is smooth if and only if $j(D) \ne 0$,
and decomposes into the union of three lines in general position
if $j(D) = 0$.
\end{proof}

\subsection{Restricting $\logf$ to other cubic curves}

Here we consider the restriction of the sheaf $\logf$
to another cubic curve $E$
defined by a polynomial $g$.
From the exact sequence (\ref{eq:exact}),
we have
$$
 0 \lto \scO(-3)|_E
   \lto \scO(-1)^{\oplus 3}|_E
   \lto \scF|_E
   \lto 0.
$$
Hence we have
$$
\begin{array}{ccccccccccc}
0  &
\rightarrow  &
H^0(\scF|_E)  &
\rightarrow  &
H^1(\scO(-3)|_E)&
\stackrel{df\wedge}{\longrightarrow}  &
H^1(\scO(-1)^{\oplus 3}|_E)&
\rightarrow  &
H^0(\scF|_E)&
\rightarrow  &
0\\
&&&&\rotatebox{90}{$\cong$}&&\rotatebox{90}{$\cong$}&&&&\\
&&&& H^0(\scO(3)|_E)^*
&
&H^0(\scO(1)^{\oplus 3}|_E)^*.&&&&
\end{array}
$$
Since $H^0(\scO(3)|_E)=\Sym^3V^*/(g)$ and
$H^0(\scO(1)|_E)^3=(V^*)^3$,
the map $df \wedge$ is dual to the map induced by
$$
\begin{array}{ccc}
 (V^*)^3 & \lto & \Sym^3 V^* \\
 \rotatebox{90}{$\in$} & & \rotatebox{90}{$\in$} \\
 (F_0, F_1, F_2) & \lto &
  F_0 \partial_0 f + F_1 \partial_1 f + F_2 \partial_2 f.
\end{array}
$$
This map is injective due to de Rham--Saito's lemma,
and the image can be identified with the degree $3$ part
$J(f)_3$ of the Jacobi ideal.
Hence we have
$$
 H^0(\scF|_E) =
  \begin{cases}
   \bC & g \in J(f)_3, \\
   0 & g \not \in J(f)_3.
  \end{cases}
$$
By an explicit calculation,
we obtain the following:
\begin{proposition} \label{prop:jumping_cubic}
Let $f_t$ be the Hesse cubic in (\ref{eq:Hesse})
and put
$$
 g = \sum_{0 \le i \le j \le k \le 2} a_{ijk} z_i z_j z_k.
$$
Then the hyperplane $J(f_t)_3 \subset \Sym^3 V^*$
is given by
$$
 J(f_t)_3 = \{ g \mid a_{012} + t (a_{000} + a_{111} + a_{222}) = 0 \}.
$$
\end{proposition}

\section{Torelli theorem}

Here we prove our main result:

\begin{theorem}
Let $C$ and $C'$ be smooth cubic curves
with non-vanishing $j$-invariants.
If $\scT(-\log C)$ is isomorphic to $\scT(-\log C')$
as an $\scO_{\bP^2}$-module, then $C = C'$.
\end{theorem}

\begin{proof}
Take a homogeneous coordinate
of the dual projective plane
so that the set of jumping lines of $\scT(- \log C)$
is a Hesse cubic.
Since a smooth cubic
whose Cayleyan curve is a smooth Hesse cubic
must be a Hesse cubic,
$C$ and $C'$ are Hesse cubics.
Then Proposition \ref{prop:jumping_cubic} shows that
$C$ must coincide with $C'$.
\end{proof}

\begin{remark}
The Torelli theorem fails for cubic curves
with vanishing $j$-invariants.
Indeed,
the family
$$
 a z_0^3 + b z_1^3 + c z_0^3 = 0, \qquad a, b, c \in \bCx
$$
consists of cubic curves
with identical Cayleyan curves given by
$$
 \alpha_0 \alpha_1 \alpha_2 = 0.
$$
Since the set of jumping lines determines
a unique stable bundle
if it consists of three lines in general position
by Barth \cite{Barth_MVBPP},
the sheaf of logarithmic vector fields
does not depend on $a$, $b$, and $c$.
\end{remark}

\bibliographystyle{plain}

\newcommand{\noop}[1]{}\def\cprime{$'$} \def\cprime{$'$}
  \def\cftil#1{\ifmmode\setbox7\hbox{$\accent"5E#1$}\else
  \setbox7\hbox{\accent"5E#1}\penalty 10000\relax\fi\raise 1\ht7
  \hbox{\lower1.15ex\hbox to 1\wd7{\hss\accent"7E\hss}}\penalty 10000
  \hskip-1\wd7\penalty 10000\box7} \def\cydot{\leavevmode\raise.4ex\hbox{.}}
  \def\cprime{$'$} \def\cprime{$'$} \def\cprime{$'$}

\noindent
Kazushi Ueda

Department of Mathematics,
Graduate School of Science,
Osaka University,
Machikaneyama 1-1,
Toyonaka,
Osaka,
560-0043,
Japan.

{\em e-mail address}\ : \  kazushi@cr.math.sci.osaka-u.ac.jp

\ \\

\noindent
Masahiko Yoshinaga

The Abdus Salam International Centre for Theoretical Physics,
Strada Costiera 11,
Trieste 34014,
Italy

{\em e-mail address}\ : \ myoshina@ictp.it

\end{document}